\documentclass[10pt]{amsart}
\usepackage{amsfonts}
\usepackage{amsmath}
\usepackage{amsthm}
\usepackage{amssymb}
\usepackage{latexsym}
\newtheorem{cor}{Corollary}[section]
\newtheorem{lem}{Lemma}[section]

\newtheorem{prop}{Proposition}[section]

\theoremstyle{definition}
\newtheorem{defn}{Definition}[section]
\theoremstyle{definition}
\newtheorem{thm}{Theorem}

\newenvironment{pf}{\proof}{\endproof}

\theoremstyle{remark}

\numberwithin{equation}{section}
\begin{document}

\newcommand{\thmref}[1]{Theorem~\ref{#1}}
\newcommand{\secref}[1]{Sect.~\ref{#1}}
\newcommand{\lemref}[1]{Lemma~\ref{#1}}
\newcommand{\propref}[1]{Proposition~\ref{#1}}
\newcommand{\corref}[1]{Corollary~\ref{#1}}
\newcommand{\remref}[1]{Remark~\ref{#1}}
\newcommand{\nc}{\newcommand}
\newcommand{\rnc}{\renewcommand}
\nc{\cal}{\mathcal}
\nc{\goth}{\mathfrak}
\rnc{\bold}{\mathbf}
\renewcommand{\frak}{\mathfrak}
\renewcommand{\Bbb}{\mathbb}

\nc{\Cal}{\mathcal}
\nc{\Xp}[1]{X^+(#1)}
\nc{\Xm}[1]{X^-(#1)}
\nc{\on}{\operatorname}
\nc{\ch}{\mbox{ch}}
\nc{\Z}{{\bold Z}}
\nc{\J}{{\cal J}}
\nc{\C}{{\bold C}}
\nc{\Q}{{\bold Q}}
\renewcommand{\P}{{\cal P}}
\nc{\N}{{\Bbb N}}
\nc\beq{\begin{equation}}
\nc\enq{\end{equation}}
\nc\lan{\langle}
\nc\ran{\rangle}
\nc\bsl{\backslash}
\nc\mto{\mapsto}
\nc\lra{\leftrightarrow}
\nc\hra{\hookrightarrow}
\nc\sm{\smallmatrix}
\nc\esm{\endsmallmatrix}
\nc\sub{\subset}
\nc\ti{\tilde}
\nc\nl{\newline}
\nc\fra{\frac}
\nc\und{\underline}
\nc\ov{\overline}
\nc\ot{\otimes}
\nc\bbq{\bar{\bq}_l}
\nc\bcc{\thickfracwithdelims[]\thickness0}
\nc\ad{\text{\rm ad}}
\nc\Ad{\text{\rm Ad}}
\nc\Hom{\text{\rm Hom}}
\nc\End{\text{\rm End}}
\nc\Ind{\text{\rm Ind}}
\nc\Res{\text{\rm Res}}
\nc\Ker{\text{\rm Ker}}
\rnc\Im{\text{Im}}
\nc\sgn{\text{\rm sgn}}
\nc\tr{\text{\rm tr}}
\nc\Tr{\text{\rm Tr}}
\nc\supp{\text{\rm supp}}
\nc\card{\text{\rm card}}
\nc\bst{{}^\bigstar\!}
\nc\he{\heartsuit}
\nc\clu{\clubsuit}
\nc\spa{\spadesuit}
\nc\di{\diamond}

\nc\al{\alpha}
\nc\bet{\beta}
\nc\ga{\gamma}
\nc\de{\delta}
\nc\ep{\epsilon}
\nc\io{\iota}
\nc\om{\omega}
\nc\si{\sigma}
\rnc\th{\theta}
\nc\ka{\kappa}
\nc\la{\lambda}
\nc\ze{\zeta}

\nc\vp{\varpi}
\nc\vt{\vartheta}
\nc\vr{\varrho}

\nc\Ga{\Gamma}
\nc\De{\Delta}
\nc\Om{\Omega}
\nc\Si{\Sigma}
\nc\Th{\Theta}
\nc\La{\Lambda}
\nc\boa{\bold a}
\nc\bob{\bold b}
\nc\boc{\bold c}
\nc\bod{\bold d}
\nc\boe{\bold e}
\nc\bof{\bold f}
\nc\bog{\bold g}
\nc\boh{\bold h}
\nc\boi{\bold i}
\nc\boj{\bold j}
\nc\bok{\bold k}
\nc\bol{\bold l}
\nc\bom{\bold m}
\nc\bon{\bold n}
\nc\boo{\bold o}
\nc\bop{\bold p}
\nc\boq{\bold q}
\nc\bor{\bold r}
\nc\bos{\bold s}
\nc\bou{\bold u}
\nc\bov{\bold v}
\nc\bow{\bold w}
\nc\boz{\bold z}

\nc\ba{\bold A}
\nc\bb{\bold B}
\nc\bc{\bold C}
\nc\bd{\bold D}
\nc\be{\bold E}
\nc\bg{\bold G}
\nc\bh{\bold h}
\nc\bH{\bold H}

\nc\bi{\bold I}
\nc\bj{\bold J}
\nc\bk{\bold K}
\nc\bl{\bold L}
\nc\bm{\bold M}
\nc\bn{\bold N}
\nc\bo{\bold O}
\nc\bp{\bold P}
\nc\bq{\bold Q}
\nc\br{\bold R}
\nc\bs{\bold S}
\nc\bt{\bold T}
\nc\bu{\bold U}
\nc\bv{\bold v}
\nc\bV{\bold V}

\nc\bw{\bold W}
\nc\bz{\bold Z}
\nc\bx{\bold x}
\nc\bX{\bold X}
\nc\blambda{{\mbox{\boldmath $\Lambda$}}}
\nc\bpi{{\mbox{\boldmath $\pi$}}}

\nc\e[1]{E_{#1}}
\nc\ei[1]{E_{\delta - \alpha_{#1}}}
\nc\esi[1]{E_{s \delta - \alpha_{#1}}}
\nc\eri[1]{E_{r \delta - \alpha_{#1}}}
\nc\ed[2][]{E_{#1 \delta,#2}}
\nc\ekd[1]{E_{k \delta,#1}}
\nc\emd[1]{E_{m \delta,#1}}
\nc\erd[1]{E_{r \delta,#1}}

\nc\ef[1]{F_{#1}}
\nc\efi[1]{F_{\delta - \alpha_{#1}}}
\nc\efsi[1]{F_{s \delta - \alpha_{#1}}}
\nc\efri[1]{F_{r \delta - \alpha_{#1}}}
\nc\efd[2][]{F_{#1 \delta,#2}}
\nc\efkd[1]{F_{k \delta,#1}}
\nc\efmd[1]{F_{m \delta,#1}}
\nc\efrd[1]{F_{r \delta,#1}}
\nc{\ug}{\bu^{fin}}

\nc\fa{\frak a}
\nc\fb{\frak b}
\nc\fc{\frak c}
\nc\fd{\frak d}
\nc\fe{\frak e}
\nc\ff{\frak f}
\nc\fg{\frak g}
\nc\fh{\frak h}
\nc\fj{\frak j}
\nc\fk{\frak k}
\nc\fl{\frak l}
\nc\fm{\frak m}
\nc\fn{\frak n}
\nc\fo{\frak o}
\nc\fp{\frak p}
\nc\fq{\frak q}
\nc\fr{\frak r}
\nc\fs{\frak s}
\nc\ft{\frak t}
\nc\fu{\frak u}
\nc\fv{\frak v}
\nc\fz{\frak z}
\nc\fx{\frak x}
\nc\fy{\frak y}

\nc\fA{\frak A}
\nc\fB{\frak B}
\nc\fC{\frak C}
\nc\fD{\frak D}
\nc\fE{\frak E}
\nc\fF{\frak F}
\nc\fG{\frak G}
\nc\fH{\frak H}
\nc\fJ{\frak J}
\nc\fK{\frak K}
\nc\fL{\frak L}
\nc\fM{\frak M}
\nc\fN{\frak N}
\nc\fO{\frak O}
\nc\fP{\frak P}
\nc\fQ{\frak Q}
\nc\fR{\frak R}
\nc\fS{\frak S}
\nc\fT{\frak T}
\nc\fU{\frak U}
\nc\fV{\frak V}
\nc\fZ{\frak Z}
\nc\fX{\frak X}
\nc\fY{\frak Y}
\nc\tfi{\ti{\Phi}}
\nc\bF{\bold F}

\nc\ua{\bold U_\A}

\nc\qinti[1]{[#1]_i}
\nc\q[1]{[#1]_q}
\nc\xpm[2]{E_{#2 \delta \pm \alpha_#1}}  
\nc\xmp[2]{E_{#2 \delta \mp \alpha_#1}}
\nc\xp[2]{E_{#2 \delta + \alpha_{#1}}}
\nc\xm[2]{E_{#2 \delta - \alpha_{#1}}}
\nc\hik{\ed{k}{i}}
\nc\hjl{\ed{l}{j}}
\nc\qcoeff[3]{\left[ \begin{smallmatrix} {#1}& \\ {#2}& \end{smallmatrix}
\negthickspace \right]_{#3}}
\nc\qi{q}
\nc\qj{q}

\nc\ufdm{{_\ca\bu}_{\rm fd}^{\le 0}}


\nc\isom{\cong} 

\nc{\pone}{{\Bbb C}{\Bbb P}^1}
\nc{\pa}{\partial}
\def\H{\cal H}
\def\L{\cal L}
\nc{\F}{{\cal F}}
\nc{\Sym}{{\goth S}}
\nc{\A}{{\cal A}}
\nc{\arr}{\rightarrow}
\nc{\larr}{\longrightarrow}

\nc{\ri}{\rangle}
\nc{\lef}{\langle}
\nc{\W}{{\cal W}}
\nc{\uqatwoatone}{{U_{q,1}}(\su)}
\nc{\uqtwo}{U_q(\goth{sl}_2)}
\nc{\dij}{\delta_{ij}}
\nc{\divei}{E_{\alpha_i}^{(n)}}
\nc{\divfi}{F_{\alpha_i}^{(n)}}
\nc{\Lzero}{\Lambda_0}
\nc{\Lone}{\Lambda_1}
\nc{\ve}{\varepsilon}
\nc{\phioneminusi}{\Phi^{(1-i,i)}}
\nc{\phioneminusistar}{\Phi^{* (1-i,i)}}
\nc{\phii}{\Phi^{(i,1-i)}}
\nc{\Li}{\Lambda_i}
\nc{\Loneminusi}{\Lambda_{1-i}}
\nc{\vtimesz}{v_\ve \otimes z^m}

\nc{\asltwo}{\widehat{\goth{sl}_2}}
\nc\eh{\frak h^e}  
\nc\loopg{L(\frak g)}  
\nc\eloopg{L^e(\frak g)} 
\nc\ebu{\bu^e} 
\nc\loopa{L(\frak a)}  

\nc\teb{\tilde E_\boc}
\nc\tebp{\tilde E_{\boc'}}

\title{On the Fermionic Formula and the 
Kirillov-Reshetikhin conjecture.}\author{Vyjayanthi Chari}
\address{Vyjayanthi Chari, Department of Mathematics, University of California, 
Riverside, CA 92521.}
\maketitle

\section{Introduction} The  irreducible finite-dimensional representations of 
quantum affine algebras $\bu_q(\hat{\frak g})$  have been studied from various 
viewpoints, \cite{AK}, \cite{CP1}, \cite{CP3}, \cite{C}, \cite{CP4}, \cite{FR}, 
\cite{FM}, \cite{KR}, \cite{K}. 
These representations decompose as  a direct sum of irreducible representations 
of the  quantized eneveloping algebra $\bu_q(\frak g)$ associated to the 
underlying finite-dimensional simple Lie algebra $\frak g$. But, except in a 
few special cases, little is known about the  isotypical components occuring in 
the decomposition. However, for  a certain class of modules (namely the one 
associated in a canonical way to a multiple of a fundamental weight of $\frak 
g$), there is a conjecture due to Kirillov and Reshetikhin \cite{KR} for 
Yangians that describes the $\frak g$-isotypical components. A combinatorial 
interpretation of their conjecture  was given by Kleber, \cite{Kl} (see also \cite{HKOTY}). It is the purpose of this paper to prove the conjecture for   
the quantum affine algebras associated to the classical simple Lie algebras, using 
Kleber's interpretation.

We now describe the conjecture and the results more explicitly. Let $\lambda_1,\lambda_2,\cdots ,\lambda_n$ be a  set of fundamental weights for $\frak g$ and, for any dominant integral weight $\mu$, let $V_q(\mu)$ denote the irreducible $\bu_q(\frak g)$-module with highest weight $\mu$. For each $m\in\bz^+$ and $i=1,\cdots ,n$, the conjecture predicts the  existence of   an irreducible representation $V^{\text{aff}}_q(m\lambda_i)$  of $\bu_q(\hat{\frak g})$ whose
highest weight when viewed as a representation of $\bu_q(\frak g)$ is
$m\lambda_i$. The decomposition of the tensor product of $N$ such 
representations as $\bu_q(\frak g)$--modules is given by,
$$
\bigotimes_{a=1}^N V_q^{\text{aff}}(m_a\lambda_{i_a})
\simeq \sum_{\lambda} n_\lambda ,V_q(\lambda)
$$
where, the sum runs over all dominant integral weights $\lambda\le \sum m_a\lambda_{i_a}$.
 The nonnegative integer
$n_\lambda$ is the multiplicity with which the irreducible $\bu_q(\frak g)$-module
$V_q(\lambda)$ occurs in the decomposition.  Write $\lambda =  \sum
m_a\lambda_{i_a} - \sum n_i \alpha_i$, $n_i\in\bz^+$.  Then
$$
n_\lambda = 
\sum_{\mbox{partitions}} \;\; \prod_{n\geq1} \;\; \prod_{k=1}^r \;\;
\left(\!\!\begin{array}{cc}P^{(k)}_n(\nu) + \nu^{(k)}_n \\ \nu^{(k)}_n
       \end{array}\!\!\right)
$$
The sum is taken over all ways of choosing partitions
$\nu^{(1)},\ldots,\nu^{(r)}$ such that $\nu^{(i)}$ is a partition of
$n_i$ which has $\nu^{(i)}_n$ parts of size $n$ (so $n_i =
\sum_{n\geq1} n \nu^{(i)}_n$).  The function $P$ is defined by
\begin{eqnarray*}
P^{(k)}_n(\nu) &=& \sum_{a=1}^N \min(n,m_a)\delta_{k,\l_a}
  - 2 \sum_{h\geq 1} \min(n,h)\nu^{(k)}_{h} + \\
&&\hspace{1cm} +
  \sum_{j\neq k}^r \sum_{h\geq 1} \min(-a_{k,j}n,-c=a_{j,k}h)\nu^{(j)}_{h}
\nonumber
\end{eqnarray*}
where $A=(a_{i,j})$ is the Cartan matrix of $\frak g$, and $\binom{a}{b}=0$
whenever $a<b$.
\vskip 6pt 
The formula describing the $n_\lambda$  is called the fermionic formula, the connection with representation theory was made by
Kirillov and Reshetikhin.  They  outlined a proof  (using the techniques of the Bethe ansatz) of the conjecture when  $\frak g$ is of type $A_n$, and showed that the  module $V_q^{\text{aff}}(m\lambda_i)$ must be isomorphic as an $A_n$--module to $V_q(m\lambda_i)$. A rigorous mathematical proof was given recently in \cite{KSS}. 

For other simple Lie algebras, the conjecture remained open, one reason being that the fermionic formula is not very tractable computationally, even in very simple cases. Although candidates were known for the modules in the case $N=m=1$,
\cite{CP3},  it was impossible to verify  the conjectures.  Kirillov and Reshetikhin did  conjecture (when $N=1$) a more explicit description of the multiplicities given by the fermionic formula.  
For instance,  when $\frak g $ is an even orthogonal algebra, and $\lambda_i$ does not correspond to the spin nodes,  then they conjectured that the   multiplicity   of $V_q(\lambda)$ in $V_q^{\text{aff}}(m\lambda_i)$ satisfies $n_\lambda \le 1$ and   \begin{equation}\label{kl} n_\lambda
\ne 0\ \ \ \text{iff}\ \ 
 \ \ \lambda=\sum_{j\ge 0}k_{i-2j}\lambda_{i-2j}, \ \ \sum_j k_{i-2j} =
m,\ \ k_r\ge 0,\end{equation}
(we understand that $\lambda_r =0$ if $r\le 0$).  This equivalence was established by Kleber \cite{Kl} who developed an algorithm to 
study the combinatorics of the fermionic formula for an arbitrary simple Lie algebra $\frak g$.  Based on this algorithm, Kleber gave  a description similar to the one above for the odd orthogonal and the symplectic Lie algebras. The exceptional cases were considered in \cite{HKOTY} where they give formulas for the multiplicities for most nodes of the Dynkin diagram.   It follows also from their work that the case of $N=1$ is the crucial case, for they prove that this implies a weak fermionic formula, which they conjecture is equivalent to the fermionic formula.

Given this explicit description of the mulitplicities,  it follows from the work of \cite{C}, \cite{CP4} on minimal affinizations
that, for any simple Lie algebra $\frak g$,  there exists up  to $\bu_q(\frak g)$--module isomorphisms, exactly one module $W_q^{\text{aff}}(m\lambda_i)=\oplus m_\mu V_q(\mu)$ which can have the prescribed decomposition.  
This is the unique minimal affinization of  $m\lambda_i$, which is characterized by the property:{\it{
$m_\mu=0$ if $m\lambda_i-\mu$ is a non--negative linear combination of simple roots which lie in a Dynkin subdiagram of type $A$}}. Thus, 
 we need to understand the $\bu_q(\frak g)$-decomposition of the 
minimal affinzations of $m\lambda_i$. We approach this problem as follows. 

In \cite{CP5}, we showed that under natural conditions, the irreducible finite-dimensional 
representations of $\bu_q(\hat{\frak g})$ admit an integral form. This allows 
us to define the $q\to 1$ limit of these representations; these are 
finite-dimensional but generally {\it reducible} representations of the loop 
algebra of 
$\frak g$. It follows by standard results that the decomposition of these 
representations of the loop algebra into a direct sum of irreducible representations of $\frak g$ is the same as the decomposition in the quantum case. In section 1, we  study the 
classical limit of the minimal affinizations and   show that for a classical simple Lie algebra $m_\mu\le 1$ and that $m_\mu\ne 0$ implies that $m_\mu$ is given by the fermionic formula.  In section 2, we  work entirely in the  quantum algebra to prove that $m_\mu=1$ 
if $\mu$  is as given in (\ref{kl}). For this, we use a result proved in \cite{K}, \cite{VV}
which describes when  a tensor product of fundamental representations of 
$\bu_q(\hat{\frak g})$ is cyclic.

Our methods also show the following for {\it any} finite-dimensional simple Lie algebra: {\it if a simple root 
$\alpha_i$ occurs with multiplicity one in the highest root of $\frak g$, then 
the modules $V_q^{fin}(m\lambda_i)$ admit a structure of a $\bu_q(\hat{\frak 
g})$-module}.  This was stated by Drinfeld in his work on Yangians, \cite{Dr1}. 
We also can prove a generalization: {\it if a root $\alpha_i$ occurs with 
multiplicity 2 in the highest root, then the minimal affinization is 
multiplicity free as a $\bu_q(\frak g)$-module}. In section 3, we summarize the results that our techniques prove for the exceptional algebras. 
\vskip 12pt
\noindent {\it Acknowledgements}. It is a pleasure to thank Michael Kleber for 
explaining his results to me. I am also grateful to M. Okado for many helpful discussions.

\section {The classical case} In this section, we study certain 
finite-dimensional modules for the loop algebra of $\frak g$. These modules (see the discussion 
following Definition \ref{wim} for their definition)  are the $q\to 1$ limit of 
 irreducible representations of the quantum loop algebra, although this does 
not become clear until the conjecture of Kirillov and Reshetikhin is 
established.  The main result of this section is Theorem \ref{main1}.

Let  $\frak{g}$ be a 
finite-dimensional complex simple Lie algebra of type $X_n$ (where $X=A,B, C$ or $D$), let $\frak h$  be a 
Cartan subalgebra of $\frak g$ and  $R$  the set of roots of $\frak g$ with 
respect to $\frak h$. Let $I=\{1,2,\cdots ,n\}$, fix a set of simple 
roots (resp. coroots) $\alpha_i$ (resp. $h_i$) ($i\in I$),  and let $R^+\subset 
\frak h^*$ be the corresponding set of positive 
roots.  We assume that the simple roots are numbered as in \cite{Bo}; in 
particular, the subset $\{j,j+1,\cdots, n\}\subset I$ defines a subalgebra of 
type $X_{n-j+1}$.

Let $Q=\bigoplus_{i=1}^n \bz\alpha_i$ (resp. 
$Q^+=\bigoplus_{i=1}^n \bz^+\alpha_i$) denote the root (resp. positive root) 
lattice 
of $\frak g$.  For $\eta\in Q^+$, $\eta=\sum_ir_i\alpha_i$, we set 
${\text{ht}}\,\eta=\sum_i r_i$.
 Let  $P$ (resp. $P^+$) be the lattice of integral (resp. dominant 
integral) weights. 
For $i\in I$, let $\lambda_i\in P^+$ be the $i^{th}$ fundamental weight. Given $\mu =\sum_{r=1}^n k_r\lambda_r\in P^+$, set $\ell(\mu) =\sum_{r=1}^n 
k_r.$

\begin{defn} For $ i\in I$ and $m\in\bz^+$ define subsets $P(i,m)$ of $P^+$ as follows:
\begin{enumerate}
\item[(i)] If $\frak g$ is of type $A_n$ then  $P(i,m) =\{m\lambda_i\}$ for all $i\in I$ and $m\in\bz^+$.
\item[(ii)] If $\frak g$ is of type $B_n$, then

\begin{align*}& P(i,1)=\{\lambda_i,\lambda_{i-2},\cdots \lambda_0\},\ \  \ \  1\le i<n,\\
&P(n,1)=\{\lambda_n\}, \ \ \ \ \ P(n,2)= \{2\lambda_n,\lambda_{n-2},\lambda_{n-4}, \cdots ,\lambda_0\},\\
&P(i,m)=P(i,1)+P(i,m-1),\ \ 1\le i<n, \ \ \ \ P(n,m)= P(n,m-2)+P(n,2),\ \ m\ge 3
\end{align*}
where $\lambda_0 =0$ if $i\in I$ is even and $\lambda_0=\lambda_1$ if $i\in I$ is odd.
\item[(iii)] If $\frak g$ is of type $D_n$, then
set
\begin{align*} P(i,1)&=\{\lambda_i,\lambda_{i-2},\cdots \lambda_0\},\ \ 1\le i<n-1\\
P(i,m)& =P(i,1)+P(i,m-1),\ \ 1\le i <n-1,n,\end{align*}
where $\lambda_0 =0$ if $i\in I$ is even and $\lambda_0=\lambda_1$ if $i\in I$ is odd. Set 
\begin{equation*} P(i,m) =\{m\lambda_i\}, \ \ i=n-1,n.\end{equation*}
\item[(iv)] If $\frak g$ is of type $C_n$, then,
\begin{align*} &P(i,1)=\lambda_i,\ \  P(i,2)=\{2\lambda_i,2\lambda_{i-1}, \cdots, 2\lambda_1, 0\},\ \ 1\le i<n,\\
&P(i,m)=P(i, m-2)+P(i,2), \ \ m\ge 3, \ \ 1\le i<n,\\
&P(n,m)= \{m\lambda_n\}.\end{align*}
\end{enumerate}
\end{defn}

\noindent The following lemma is trivially checked.
\begin{lem}
\begin{enumerate}\item[(i)] If $\frak g$ is of type $B_n$ and $1\le i<n$, then,
\begin{align*} P(i,m) &=\left\{\sum_{j=0}^{[i/2]} k_{i-2j}\lambda_{i-2j}: \sum_j k_{i-2j} =m\right\},\\
P(n,m)&=\left\{\sum_{j=0}^{[i/2]} k_{i-2j}\lambda_{i-2j}: k_n+2\sum_j k_{i-2j} =m\right\}.\end{align*}
\item[(ii)] if $\frak g$ is of type $D_n$, and $1\le i\le n-2$, then 
\begin{equation*} P(i,m) =\left\{\sum_{j=0}^{[i/2]} k_{i-2j}\lambda_{i-2j}: \sum_jk_{i-2j} =m\right\}.\end{equation*}
\item[(iii)] If $\frak g$ is of type $C_n$, then we set $\lambda_0=0$ and  
\begin{equation*} P(i,m)=\left\{\sum_{j=0}^i k_j\lambda_j: \sum_jk_j=m,\ \ k_i\equiv m {\text{mod}} \ 2, \ k_j\equiv 0 {\text{mod}}\ 2, j\ne i\right\}.\end{equation*}\end{enumerate}\hfill\qedsymbol
\end{lem}

Let  $\frak n^\pm$ be the subalgebras 
\begin{equation*} \frak n^\pm =\bigoplus_{\pm\alpha\in R^+}  \frak 
g_\alpha.\end{equation*}
{\it Throughout this paper we shall (by abuse of notation) denote any non-zero 
element of $\frak g_{\pm\alpha}$ as $x_{\alpha}^\pm$; of course, any two such 
elements are scalar multiples of each other, but for our purposes  a precise 
choice of scalars is irrelevant. Thus, if $\alpha,\beta\in R^+$ is such that 
$\alpha
\pm \beta\in R^+$, then we shall write
\begin{equation*} [x_\alpha^+,x_\beta^\pm] =x_{\alpha\pm\beta}^+,\end{equation*}
etc.}

\medspace

For any Lie algebra $\frak a$, the   loop algebra of $\frak a$ is the Lie 
algebra
\begin{equation*}\loopa = \frak a\otimes \bc[t,t^{-1}],\end{equation*}
with commutator given by
\begin{equation*}
 [x\otimes t^r, y\otimes t^s] 
=[x,y]\otimes t^{r+s},\end{equation*}
for $x,y\in\frak a$, $r,s\in\bz$. 
For any $x\in \frak a$, $m\in\bz$, we denote by $x_m$ the element $x\otimes 
t^m\in\loopa$.  Let $\bu(\frak a)$ be  the universal enveloping algebra of $\frak a$. 

For $i\in I$, $k\in\bz$, define elements of $L(\frak g)$ by   
$e_i^\pm =x_{\alpha_i}^\pm\otimes 1$, $x_{i,k}^\pm=x_{\alpha_i}^\pm \otimes t^k$ and $e_0^\pm =x_{\theta_1}^\mp\otimes t^{\pm 1}$. Then, the elements $e_i^\pm$  
($i=0,\cdots ,n$) generate $L(\frak g)$.  
 We set 
\begin{equation*}
\ \bu(\loopg)= \bu,\ \ \bu(\frak g)=\ug. 
\end{equation*}
We have
\begin{equation*}\bu=\bu(L(\frak n^-))\bu(L(\frak h))\bu(L(\frak n^+)),\ \  
\bu(\frak g)=\bu(\frak n^-)\bu(\frak h)\bu(\frak n^+).\end{equation*}

Given  $\lambda\in P^+$, let $V^{fin}(\lambda)$ be the unique irreducible 
finite-dimensional $\bu^{fin}$-module with highest weight $\lambda$ with highest weight vector  
$ v_\lambda$. For all $\alpha\in 
R^+$, $h\in \frak h$, we have 
\begin{equation*} x_\alpha^+. v_\lambda= 0, \ \ h.v =\lambda(h).v_\lambda,\ \ 
(x_\alpha^-)^{\lambda(h_\alpha)+1}.v_\lambda =0.\end{equation*}
  The action of $\frak g$ on $V^{fin}(\lambda)$ extends to an action of $\loopg$, by setting,\begin{equation*} x_m. v= x.v,\ \ \forall\ m\in\bz, \ x\in\frak g.\end{equation*}
We denote this $\loopg$--module by $V(\lambda)$.
For any finite-dimensional $\ug$-module $V$ and  any  $\nu\in \frak h^*$, let
\begin{equation*} V_\nu=\{v\in V: h.v=\nu(h)v\ \forall \ h\in\frak 
h\}.\end{equation*} 
Since $V$ is a direct sum of irreducible $\ug$-modules, we can write
\begin{equation*} V\cong \bigoplus_{\mu\in P^+}m_\mu(V) 
V^{fin}(\mu),\end{equation*}
where $m_\mu(V)\ge 0$ is the multiplicity with which $V^{fin}(\mu)$ occurs in  
the sum.

We next recall the definition of certain highest weight modules, introduced in 
\cite{CP5}; in fact, only the following special case will be needed.
Let $\bpi_{i,m}= (\pi_1,\cdots ,\pi_n)$ be the $n$-tuple of polynomials in 
$\bc[u]$ given by
\begin{equation*}  \pi_j(u)
 =1\ \ {\text{if}\ } \ j\ne i,\ \ \ 
\pi_i(u)= (1-u)^m.\end{equation*}
\begin{defn}{\label{wim}}The $\bu$-modules $W(\bpi_{i,m})$ are generated by an 
element $w_{i,m}$ subject to the relations
\begin{align}\label{relw0} x_{j,k}^+.w_{i,m} =0,\ \ & \ \  h_{k}.w_{i,m} = m\lambda_i(h) 
w_{i,m} \ \ \ (h\in\frak h,\  k\in\bz), \\
\label{relw} (x_{i,k}^-)^{m+1}.w_{i,m}= 0,\ \ &\ \ x_{j,k}^-.w_{i,m} 
=0\ \ \ \ \ (j\ne i,\  k\in\bz).\end{align}\hfill\qedsymbol\end{defn}
 
The following proposition was proved in \cite[Section2,  Theorem 1]{CP5}. 
\begin{prop}{\label{cpw}} The $\bu$-module $W({\bpi_{i,m}})$ is finite 
dimensional and
\begin{equation*} \bu(\frak n^-\otimes\bc[t]).w_{i,m} 
=W({\bpi_{i,m}}).\end{equation*}
Further, the module $V(m\lambda_i)$ is the 
unique irreducible quotient of $W(\bpi_{i,m})$. In particular, 
$m_{m\lambda_i}(W({\bpi_{i,m}})) =1$. 
\hfill\qedsymbol\end{prop}

The elements
\begin{equation}{\label{reliw}} x_{i,k}^-.w_{i,m}-x_{i,0}^-.w_{i,m} \ \ \ \ \ \ 
(k\in\bz) \end{equation} 
generate a proper $\bu$-submodule of $W(\bpi_{i,m})$. Let $W(i,m)$ denote the quotient 
of $W(\bpi_{i,m})$ by this submodule.  
 We continue to denote by  $w_{i,m}$  the image of $w_{i,m}\in 
W(\bpi_{i,m})$ in $W_{i,m}$. 
The main result of this section is the following.
\medskip

\begin{thm}{\label{main1}} Let $i\in I$, $m\ge 0$. For all $\mu\in P^+$ we have $m_{\mu}(W(i,m))\le 1$. Further, \begin{equation*} m_\mu(W(i,m))\ne 0\implies \mu\in P(i,m).\end{equation*}
\end{thm}
The rest of the section is devoted to proving the theorem.

\medskip

For $i\in I$ and $l=0,1,2$, set
\begin{equation*} R(i,l) =\left\{\sum_{k=1}^n m_k\alpha_k\in R^+ : m_i 
=l\right\}.\end{equation*}
Clearly, 
\begin{equation*} R^+ =\bigcup_{l=0}^2 R(i,l).\end{equation*}
For $i\in I$, $l=0,1,2$, define the subspaces  $\frak n^\pm(i,l)$ in the 
obvious way. Then,
\begin{align*} [\frak{n}^\pm(i,l'), \frak{n}^\pm(i,l)] &= 0,\ \  \ {\text{if}}\ 
l'+l>2,\\ 
[\frak{n}^\pm(i,l'), \frak{n}^\pm(i,l)]& =\frak{n}^\pm (i,l'+l),\ \    \ 
{\text{if}}\ l'+l\le 2.\end{align*}

\begin{prop}\label{limits} Let $\alpha\in R^+$, $f\in\bc[t,t^{-1}]$. Then, 
\begin{equation*}\alpha\in R(i,l)\implies (x_\alpha^-\otimes f(t-1)^l) .w_{i,m} 
=0.\end{equation*}
\end{prop}
\begin{pf} We proceed by induction on ${\text{ht}}\,\alpha$. The case of 
${\text{ht}}\,\alpha=1$ is  clear from (\ref{relw}) and (\ref{reliw}).
Assume that the result holds for   ${\text{ht}}\,\alpha <r$.  Choosing  $j\in I$ so that $\beta=\alpha-\alpha_j\in R^+$, we get
\begin{equation*} x^-_\alpha\otimes fg =[x^-_{\alpha_j}\otimes f, x^-_{\beta}\otimes g],\end{equation*}
 for all $f,g\in\bc[t,t^{-1}]$.

If $j\ne i$, then $\alpha,\beta\in R(i,l)$ for some $l=0,1, 2$. 
Now  (\ref{relw}) gives \begin{equation*} (x_{\alpha}^-\otimes f(t-a)^l).w_{i,m} = 
(x^-_{\alpha_j}x^-_{\beta}\otimes f(t-a)^l).w_{i,m}.\end{equation*}
Since ${\text{ht}}\,\beta< {\text{ht}}\,\alpha$, the result follows.
Assume now that $j=i$.   
If $\alpha\in R(i,1)$, then $\beta\in R(i,0)$ and we get by using  induction and \eqref{reliw} 
that \begin{equation*}(x_{\alpha}^-\otimes f(t-1)).w_{i,m} = 
-(x_{\beta}^-.x_{\alpha_i}^-\otimes f(t-1)).w_{i,m} =0.\end{equation*} 
Finally, if $\alpha\in R(i,2)$,  then $\beta\in R(i,1)$ and we have again by induction 
that    
\begin{equation*} (x_{\alpha}^-\otimes f(t-1)^2).w = [x_{\alpha_i}^-\otimes 
(t-1),  x^-_{\beta}\otimes f(t-1)].w_{i,m} =0. \end{equation*}
This proves the proposition.
\end{pf}

The following is now  immediate by applying the PBW theorem.
\begin{cor}{\label{r=1}} We have,
\begin{equation*} W(i,m)=\bu(\frak n^-)\bu(\frak n^-(i,2)\otimes (t-1)). w_{i,m}.\end{equation*} In particular if $R(i,2)=\{\phi\}$  then 
\begin{equation*} W(i,m)\cong V(m\lambda_i)\ \ \forall\  m\in\bz_+.\end{equation*}\hfill\qedsymbol
\end{cor}
In view of this corollary, we can now assume that $\frak g$ is of type $B_n$, $C_n$ or $D_n$ and that $i\ne 1$ (resp. $i\ne n$, $i\ne 1, n-1, n$). 
We list the sets $R(i,2)$ explicitly in these cases. Define roots,
\begin{align*} \theta^i_{l,k}&=\sum_{j=l}^k\alpha_j+2\sum_{j=k+1}^n\alpha_j.\ \ \text{if $\frak g= B_n$}, \ \ 1\le l\le k\le i-1,\\
&= \sum_{j=l}^k\alpha_j+2\left(\sum_{j=k+1}^{n-2}\alpha_j \right)+\alpha_{n-1}+\alpha_n,\ \ \text{if $\frak g =D_n$},\ \ 1\le l\le k\le i-1,\\
&= \sum_{j=l}^{k-1}\alpha_j+2\left(\sum_{j=k}^{n-1}\alpha_j\right)+\alpha_n,\
\ \text{if $\frak g=C_n$},\ \ 1\le l\le k\le i.\end{align*}
The collection of all the $\theta^i_{k,l}$ is $R(i,2)$. Let $\frak u^i$ be the subalgebra of $\frak g$ spanned by $\{x^-_{\theta^i_{j,j}}: 1\le j\le i-1, \ \ i-1\equiv j\mod 2\}$ (resp. $ \{x^-_{\theta^i_{j,j}}: 1\le j\le i \}$), if $\frak g$ is of type $B_n$ or $D_n$ (resp. $C_n$).  

To prove the next  proposition only, we shall denote  by $\frak g_n$ the Lie algebra of type $X_n$ and by $W_n(i,m)$ the module $W(i,m)$ etc. The assignment 
\begin{equation*} x^\pm_{\alpha_j}\to x^\pm_{\alpha_{j+1}},\end{equation*}
extends to an embedding of $\frak g_{n-1}\to \frak g_n$ and to the corresponding loop algebras.   Let $\frak t^i_n =\oplus_{k}\bc x^-_{\theta^i_{1,k}}$ and let $\frak{n}_{n-1}(i,2)$ denote the image in $\frak{n}_n$ of $\frak{n}_{n-1}(i-1,2)$ etc. Then, \begin{equation*}
\frak n^-_n(i,2) =\frak n_{n-1}^-(i,2)\oplus \frak t^i_n\ \ \frak u_n^i=\frak u_{n-1}^i\oplus \bc x^-_{\theta^i_{1,1}}.\end{equation*}
Further, it is easy to see that there exists a $\bu_{n-1}$--module map $W_{n-1}(i-1,m)\to W_n(i,m)$, for $i\in I_n$, $i>1$ (and as stated earlier $i\ne n$ for $C_n$ and  $i\ne n-1,n$ for $D_n$) with image $\bu_{n-1}.w_{i,m}$.

 We now prove,
\begin{prop} We have,
\begin{equation*} W_n(i,m)= \bu(\frak n^-_n)\bu(\frak u_n^i\otimes (t-1)).w_{i,m}.\end{equation*}
 \end{prop}
\begin{pf} 
 
We prove this proposition by induction on $n$.
In the case when $R(i,2)$ consists of exactly one element, we have $\frak u_n^i=\frak n_n^-(i,2)$ and the result is just Corollary \ref{r=1}. Hence the proposition is established for $B_2=C_2$, for  $D_4$ and for $i=1$ for all $C_n$.

So to complete the inductive step, we can assume that
$i>1$ and that  the result holds for $\frak g_{n-1}$.
Thus the induction hypothesis gives,
\begin{equation*} \bu_{n-1}.w_{i,m} = \bu(\frak{n}_{n-1}^-)\bu(\frak u^i_{n-1}\otimes (t-1)).w_{i,m}\end{equation*}
We now get,
\begin{align*} W_n(i,m) &=\bu(\frak n_n^-)\bu(\frak t^i_n\otimes (t-1))\bu(\frak n^-_{n-1}(i,2)\otimes (t-1)).w_{i,m} \\&= \bu(\frak n_n^-)\bu(\frak t^i_n\otimes (t-1))\bu(\frak n^-_{n-1})\bu(\frak u^i_{n-1}\otimes (t-1)).w_{i,m}.\end{align*}
Since  $[\frak t^i_n, \frak n^-_n]\subset \frak t_n^i$, we get
\begin{equation*} W_n(i,m)= \bu(\frak n_n^-)\bu(\frak t^i_n\otimes (t-1))\bu(\frak u^i_{n-1}\otimes (t-1)).w_{i,m}.\end{equation*} 
To complete the proof, we must show that \begin{equation}\label{crux}
\bu(\frak t^i_n\otimes (t-1))\bu(\frak u^i_{n-1}\otimes (t-1)).w_{i,m}\subset \bu(\frak n^-)\bu(\frak u^i_n\otimes (t-1)).w_{i,m}.\end{equation}

We do this in the case of $D_n$ and when $i$ is even, the proof in the other cases, is similar and simpler. Set $\theta^i_{l,k}= \theta_{l,k}$ and define elements $\gamma_j\in R^+$ by, \begin{align*}&\gamma_j =\theta_{1,j}-\theta_{j,j} =\sum_{r=1}^{j-1}\alpha_r\ \ \ {\text{if $j$ is odd}},\\
&\gamma_j=\theta_{1,j}-\theta_{j+1j+1}=\sum_{r=1}^{j+1}\alpha_r,\ \  {\text{if $j$ is even}}.\end{align*}
Since $i$ is even, we have $x^-_{\gamma_j}.w_{i,m}=0$ for all $2\le j\le i-1$.
Now, a simple checking shows that 
\begin{align*}&(x^-_{\gamma_2})^{s_2}(x^-_{\gamma_3})^{s_3}\cdots 
(x^-_{\gamma_{i-1}})^{s_{i-1}}\\&\times (x^-_{\theta_{1,1}}\otimes (t-1))^{r_1}(x^-_{\theta_{3,3}}\otimes (t-1))^{r_2+r_3+s_3}
\cdots 
(x^-_{\theta_{i-1,i-1}}\otimes (t-1))^{r_{i-2}+r_{i-1}+s_{i-1}}.w_{i,m}\\
&=(x^-_{\theta_{1,1}}\otimes (t-1))^{r_1}[(x_{\gamma_2}^-)^{r_2}(x_{\gamma_3}^-)^{r_3}, 
(x^-_{\theta_{3,3}}\otimes (t-1))^{r_2+r_3+s_3}]\cdots\\&\times 
[(x_{\gamma_{i-2}})^{r_{i-2}}(x_{\gamma_{i-1}})^{r_{i-1}} 
,x^-_{\theta_{i-1,i-1}}\otimes (t-1))^{r_{i-2}+r_{i-1}+s_{i-1}}].w_{i,m}\\
&= (x^-_{\theta_{1,1}}\otimes (t-1))^{r_1}(x^-_{\theta_{1,2}}\otimes (t-1))^{r_2}\cdots (x^-_{\theta_{1,i-1}}\otimes (t-1))^{r_{i-1}}
\\&\times
(x^-_{\theta_{3,3}}\otimes (t-1))^{s_3}
(x^-_{\theta_{5,5}}\otimes (t-1))^{s_5}\cdots 
(x^-_{\theta_{i-1,i-1}}\otimes (t-1))^{s_{i-1}}.w_{i,m},\end{align*}
where the last equality follows from the definition of the $\gamma_j$'s  and noting that $\theta_{j,j}+\gamma_k+\gamma_l\notin R^+$. This clearly proves \eqref{crux} and the proof of the proposition is complete.

\end{pf}

\noindent{\it Proof of Theorem   \ref{main1}.} Set $l=\text{dim}\ \frak u^i$ and let $\le $ be the lexicographic ordering on $\bz_+^l$. Given $\bos\in\bz_+^l$, let \begin{equation*} \bx_\bos= \prod_{j=1}^{i-1}(x^-_{\theta_{j,j}}\otimes (t-1))^{s_j},\end{equation*}
if $\frak g$ is of type $C_n$, the coresponding analogues for $B_n$ and $D_n$ are defined in the obvious way.

Let $W_\boo$ be the $\frak g$--submodule of $W(i,m)$ generated by $w_{i,m}$ and let $W_1$ be a $\frak g$--module such that
\begin{equation*} W(i,m)= W_\boo\oplus W_1.\end{equation*} 
If $W_1\ne 0$,  choose $\bos_1$ minimal so that the element $\bx_{\bos_1}.w_{i,m}$ has a non--zero projection $w_{\bos_1}$ onto $W_1$. Now choose a $\frak g$--submodule $W_2$ of $W_1$ so that,
\begin{equation*} W_1= \bu(\frak g).w_{\bos_1}\oplus W_2.\end{equation*}
Repeating, we see that we can find a finite number of elements, say  $\{w_{\bos_j}:1\le j\le k\}$, with $\bos_1<\bos_2<\cdots <\bos_k$ such that 
\begin{equation*} W(i,m) = W_\boo\oplus W_{\bos_1}\oplus\cdots\oplus W_{\bos_k},\end{equation*}
where $W_{\bos_j}=\bu(\frak g).w_{\bos_j}$.
Notice that by choice,  the projection of $\bx_{\bos}.w_{i,m}$ onto $W_{\bos_j}$ is zero if $\bos<\bos_j$.  We claim that,
\begin{equation}{\label{hw}}x_{\alpha}^+.w_{\bos_j} =0,\ \ \forall\ \  
\alpha\in R^+. \end{equation}

From now on, we assume that $\frak g$ is of type $C_n$,  the proof in the other cases is similar.
Thus, notice that if $k\ne i$, we have 
\begin{align*} x_{\alpha_k}^+.\bx_{\bos_j}.w_{i,m}&= 0\ \ \text{if}\ \ k>i,\\
&= 
x^-_{\theta_{k,k}-\alpha_k}\otimes (t-1)(x^-_{\theta_{k,k}}\otimes (t-1))^{s_k-1} \prod_{j'\ne k}(x^-_{\theta_{j',j'}}\otimes (t-1))^{s_{j'}}.w_{i,m}, \\
&= x^-_{\alpha_k}(x^-_{\theta_{k,k}}\otimes (t-1))^{s_k-1}(x^-_{\theta_{k+1,k+1}}\otimes (t-1))^{s_{k+1}+1} \prod_{j'\ne k, k+1}(x^-_{\theta_{j',j'}}\otimes (t-1))^{s_{j'}}.w_{i,m}.\end{align*}
But the right hand side of the last equality is clearly in $\oplus W_{\bos_r}$ with $\bos_r<\bos_j$. This gives \eqref{hw} if $k\ne i$. If $k=i$, then,
we have,
\begin{align*} x_{\alpha_i}^+.\bx_{\bos_j}.w_{i,m}&= x^-_{\theta_{i,i}-\alpha_i}\otimes (t-1)(x^-_{\theta_{k,k}}\otimes (t-1))^{s_k-1} \prod_{j'\ne k}(x^-_{\theta_{j',j'}}\otimes (t-1))^{s_{j'}}.w_{i,m},\\
&= 0,
\end{align*}
where the last equality follows from the fact that $\theta_{i,i}-\alpha_i\in R(i,1)$. This proves \eqref{hw} completely  and 
 and hence we get that 
if $ m_\mu(W(i,m))\ne 0$ then $\mu =m\lambda_i$ or $\mu$ is the weight of the element $w_{\bos_j}$ for some $j$. 

A simple calculation  shows that
$\theta_{j,j} =2\lambda_j-2\lambda_{j-1}$ and hence the weight of the element $w_{\bos_j}$ where $\bos_j =(s_{j1}, s_{j2},\cdots s_{jl})$ is \begin{equation*}\mu_j=
(m-2s_{ji})\lambda_i + 2(s_{ji}-s_{ji-1})\lambda_{i-1}+\cdots +2(s_{j2}-s_{j1})\lambda_1.\end{equation*} Since $\mu_j$ 
must be a dominant integral weight we see that $\mu_j\in P(i,m)$. Further, the $\mu_j$ are clearly distinct and hence Theorem \ref{main1} is proved.

\hfill\qedsymbol

 \section{The quantum case} In this section we recall the definition of the 
quantum affine algebras and several results on the irreducible 
finite-dimensional representations of $\bu_q(\hat{\frak g})$. We then define 
the module whose decomposition we are interested in and establish the 
Kirillov-Reshetikhin conjecture in this case. We continue to assume that $\frak g$ is of type $X_n$, where $X =A,B, C$ or $D$.

Let $q$ be an indeterminate, let $\bc(q)$ be the field of rational
functions in $q$ with complex coefficients, and let $\ba=\bc[q,q^{-1}]$ be
the subring of Laurent polynomials. For $r,m\in\bn$, $m\ge r$, define
\begin{equation*} 
[m]=\frac{q^m -q^{-m}}{q -q^{-1}},\ \ \ \ [m]! =[m][m-1]\ldots [2][1],\ \ \ \ 
\left[\begin{matrix} m\\ r\end{matrix}\right] 
= \frac{[m]!}{[r]![m-r]!}.
\end{equation*}
Then, $\left[\begin{matrix} m\\r\end{matrix}\right]\in\ba$. 

We now recall the definition of the quantum affine algebra.  Let $\hat A 
=(a_{ij})$ be  the $(n+1)\times (n+1)$ extended Cartan matrix associated to 
$\frak g$. Let $\hat{I} =I\cup\{0\}$. Fix non--negative integers $d_i$, $i\in\hat{I}$ such that the matix $(d_ia_{ij})$ is symmetric. Set $q_i=q^{d_i}$ and  $[m]_i=[m]_{q_i}$.
\begin{prop}{\label{defnbu}} There is a Hopf algebra $\tilde{\bu}_q$ over 
$\bq(q)$ which is generated as an algebra by elements $E_{\alpha_i}$, 
$F_{\alpha_i}$, $K_i^{{}\pm 1}$ ($i\in\hat I$), with the following defining 
relations:
\begin{align*} 
  K_iK_i^{-1}=K_i^{-1}K_i&=1,\ \ \ \ K_iK_j=K_jK_i,\\ 
  K_iE_{\alpha_j} K_i^{-1}&=q_i^{ a_{ij}}E_{\alpha_j},\\ 
K_iF_{\alpha_j} K_i^{-1}&=q_i^{-a_{ij}}F_{\alpha_j},\\
  [E_{\alpha_i}, F_{\alpha_j}
]&=\delta_{ij}\frac{K_i-K_i^{-1}}{q_i-q_i^{-1}},\\ 
  \sum_{r=0}^{1-a_{ij}}(-1)^r\left[\begin{matrix} 1-a_{ij}\\ 
  r\end{matrix}\right]_i
&(E_{\alpha_i})^rE_{\alpha_j}(E_{\alpha_i})^{1-a_{ij}-r}=0\ 
  \ \ \ \ \text{if $i\ne j$},\\
\sum_{r=0}^{1-a_{ij}}(-1)^r\left[\begin{matrix} 1-a_{ij}\\ 
  r\end{matrix}\right]_i
&(F_{\alpha_i})^rF_{\alpha_j}(F_{\alpha_i})^{1-a_{ij}-r}=0\ 
  \ \ \ \ \text{if $i\ne j$}.
\end{align*}
The comultiplication of ${\tilde{\bu}_q}$ is given on generators by
$$\Delta(E_{\alpha_i})=E_{\alpha_i}\ot 1+K_i\ot E_{\alpha_i},\ \ 
\Delta(F_{\alpha_i})=F_{\alpha_i}\ot K_i^{-1} + 1\ot F_{\alpha_i},\ \ 
\Delta(K_i)=K_i\ot K_i,$$
for $i\in\hat I$.\hfill\qedsymbol
\end{prop}
Set $K_{\theta} =\prod_{i=1}^n K_i^{r_i/d_i}$, where $\theta=\sum 
r_i\alpha_i$ is the highest root in $R^+$. 
Let $\bu_q$ be the  quotient of $\tilde{\bu}_q$ by the ideal generated by the 
central element  $K_0K_{\theta}^{-1}$; we call this the quantum loop algebra 
of $\frak g$.

It follows from 
\cite{Dr2}, \cite{B}, \cite{J} that $\bu_q$ is isomorphic to the 
algebra with generators $\bx_{i,r}^{{}\pm{}}$ ($i\in I$, $r\in\bz$), $K_i^{{}\pm 
1}$ 
($i\in I$), $\bh_{i,r}$ ($i\in I$, $r\in \bz\backslash\{0\}$) and the following 
defining relations:
\begin{align*}
   K_iK_i^{-1} = K_i^{-1}K_i& =1, \ \  
 K_iK_j =K_jK_i,\\  K_i\bh_{j,r}& =\bh_{j,r}K_i,\\ 
 K_i\bx_{j,r}^\pm K_i^{-1} &= q_i^{{}\pm
    a_{ij}}\bx_{j,r}^{{}\pm{}},\ \ \\ 
  [\bh_{i,r},\bh_{j,s}]=0,\; \; & [\bh_{i,r} , \bx_{j,s}^{{}\pm{}}] =
  \pm\frac1r[ra_{ij}]_i\bx_{j,r+s}^{{}\pm{}},\\ 
 \bx_{i,r+1}^{{}\pm{}}\bx_{j,s}^{{}\pm{}} -q_i^{{}\pm
    a_{ij}}\bx_{j,s}^{{}\pm{}}\bx_{i,r+1}^{{}\pm{}} &=q_i^{{}\pm
    a_{ij}}\bx_{i,r}^{{}\pm{}}\bx_{j,s+1}^{{}\pm{}}
  -\bx_{j,s+1}^{{}\pm{}}\bx_{i,r}^{{}\pm{}},\\ [\bx_{i,r}^+ ,
  \bx_{j,s}^-]=\delta_{i,j} & \frac{ \psi_{i,r+s}^+ -
    \psi_{i,r+s}^-}{q_i - q_i^{-1}},\\ 
\sum_{\pi\in\Sigma_m}\sum_{k=0}^m(-1)^k\left[\begin{matrix}m\\k\end{matrix}
\right]_i
  \bx_{i, r_{\pi(1)}}^{{}\pm{}}\ldots \bx_{i,r_{\pi(k)}}^{{}\pm{}} &
  \bx_{j,s}^{{}\pm{}} \bx_{i, r_{\pi(k+1)}}^{{}\pm{}}\ldots
  \bx_{i,r_{\pi(m)}}^{{}\pm{}} =0,\ \ \text{if $i\ne j$},
\end{align*}
for all sequences of integers $r_1,\ldots, r_m$, where $m =1-a_{ij}$, $\Sigma_m$ is the symmetric group on $m$ letters, and the $\psi_{i,r}^{{}\pm{}}$ are 
determined by equating powers of $u$ in the formal power series 
$$\sum_{r=0}^{\infty}\psi_{i,\pm r}^{{}\pm{}}u^{{}\pm r} = K_i^{{}\pm 1} 
{\text{exp}}\left(\pm(q_i-q_i^{-1})\sum_{s=1}^{\infty}\bh_{i,\pm s} u^{{}\pm 
s}\right).$$ 
\vskip 12pt

For $i\in I$, the above isomorphism maps $E_{\alpha_i}$ to  $\bx_{i,0}^+$ and 
$F_{\alpha_i}$ to $\bx_{i,0}^-$.  The subalgebra generated by $E_{\alpha_i}$, $F_{\alpha_i}$, $i\in I$, is the quantized enveloping algebra $\bu_q^{fin}$ associated to $\frak g$,

Define the $q$-divided powers
\begin{equation*}(\bx_{i,k}^\pm)^{(r)} 
=\frac{(\bx_{i,k}^\pm)^r}{[r]_i!},\end{equation*}
for all $i\in I$, $k\in\bz$, $r\ge 0$.
The elements $E_{\alpha_i}^{(r)}$ etc. are defined similarly. Let $\bu_\ba$ be the $\ba$-subalgebra of $\bu_q$ generated by the $K_i^{\pm 
1}$, $(\bx_{i,k}^\pm)^{(r)}$  ($i\in I$, $k\in\bz$,  $r\ge 0$).

 \begin{lem} The subalgebra $\bu_\ba$ is an $\ba$--lattice in $\bu_q$, and \begin{equation*}
\bu_q=\bc(q)\otimes_\ba\bu_\ba.
\end{equation*}
\end{lem}
\begin{pf} Let $\tilde\bu_\ba$ be the $\ba$--subalgebra generated by the elements $E_{\alpha_i}^{(r)}$, $F_{\alpha_i}^{(r)}$, $i\in\hat{I}$. It is proved in \cite{L3} that $\tilde\bu_\ba$ is an $\ba$--lattice and that
\begin{equation*}
\bu_q=\bc(q)\otimes_\ba\tilde\bu_\ba.
\end{equation*}
Hence to prove the lemma it suffices to show that the $\bu_\ba=\tilde\bu_\ba$.
For this, in view of the isomorphism between the two presentations it suffices to show that the elements $E_{\alpha_0}^{(r)}$ and $F_{\alpha_0}^{(r)}$ are in $\bu_\ba$. In the simply laced case this was proved in \cite[Proposition 2.6]{BCP}. The proof given there works as long as there  exists a simple root $\alpha_{i_0}$ which occurs with mulitplicity one in $\theta$, i.e $r_{i_0}=1$. An inspection shows that this is true  for the classical simple Lie algebras.\end{pf} 

Given $i,j\in I$ with $a_{ij}= -2$ and $k,l\in\bz$, it is easy to see that the subalgebra generated by the elements $\bx_{i,k}^\pm$ and $\bx_{j,l}^\pm$ is isomorphic to the quantized enveloping algebra of  $\bu_q(sp_5)$. Define elements,
\begin{equation*}\gamma_{k,l}(q) = \bx_{i,k}^-\bx_{j,l}^- -q^2\bx_{j,l}^-\bx_{i,k}^-,\ \ (\gamma_{k,l}(q))^{(r)}=\frac{(\gamma_{k,l}(q))^r}{[r]_i!},\end{equation*}
and 
\begin{equation*}\gamma'_{i,k}(q)=[\bx_{i,l}^-, \gamma_{k,l}(q)],\ \ (\gamma'_{k,l}(q))^{(r)} = \frac{(\gamma'_{k,l}(q))^r}{[r]_j!}.\end{equation*}
It is easy to see using the defining relations in $\bu_q$ that,\begin{equation*}\label{aij}\gamma_{k,l}(q)=q^2\gamma_{k-1,l+1}(q^{-1}),\ \ \gamma'_{k,l}(q)=q^2\gamma'_{k-1,l+1}(q^{-1}).\end{equation*}
\begin{lem}\label{aij1} Assume that $i,j\in I$ is such that $a_{ij}=-2$. Then, 
 \begin{equation*} (\bx_{i,k}^-)^{(a)}(\bx_{j,l}^-) ^{(b)} =\sum_{r,t\in \bz_+} f_{r,t} (\bx_{j,l}^-) ^{(b-r-t)}(\gamma_{k,l}(q))^{(r)}(\gamma'_{k,l}(q))^{(t)}
(\bx_{i,k}^-)^{(a-r-2t)},\end{equation*}
where $f_{r,t}\in q^{\bz_+}$. In particular the elements $(\gamma_{k,l}(q^{\pm 1}))^{(r)}$ and $(\gamma'_{k,l}(q^{\pm 1}))^{(r)}$ are in $\bu_\ba$.  \end{lem}
\begin{pf} This follows from the result proved in \cite{L3} for the quantized enveloping algebra of $sp_5$. \end{pf}

For any $\bu_q^{fin}$-module $V_q$ and any $\mu\in P$, set
\begin{equation*}( V_q)_\mu=\{ v\in V_q: K_i.v =q_i^{\mu(h_i)}v ,\ \ 
\forall \ i\in I\}.\end{equation*}
We say that $V_q$ is a module of type 1 if
\begin{equation*} V_q=\bigoplus_{\mu\in P}(V_q)_\mu.\end{equation*}
From now on, we shall only be working with $\bu_q^{fin}$-modules of type 1. 

The irreducible finite-dimensional $\bu_q^{fin}$-modules  are 
parametrized by $P^+$. Thus, for each $\lambda\in P^+$, there exists a unique 
irreducible finite-dimensional module $V_q^{fin}(\lambda)$ generated by a 
non-zero  element $v_\lambda$, with defining relations
\begin{equation*} \bx_{i,0}^+. v_\lambda= 0,\ \ K_i. 
v_\lambda=q^{\lambda(h_{i})}v_\lambda,\ \ 
(\bx_{i,0}^-)^{\lambda(h_i)+1}.v_\lambda=0,\ \ \forall \ i\in I.\end{equation*}
Further,  \begin{equation*} (V_q^{fin}(\lambda))_\mu\ne 0\implies 
\mu\in\lambda-Q^+.\end{equation*}
Set $V^{fin}_\ba(\lambda)= \bu_\ba.v_\lambda$. Then,
\begin{equation*}
V_q^{fin}(\lambda)=\bc(q)\otimes_\ba V^{fin}_\ba(\lambda),
\end{equation*}and 
\begin{equation*}\overline{V^{fin}_q(\lambda)} = \bc_1\otimes_\ba 
V^{fin}_\ba(\lambda).\end{equation*} Then \cite{L1}, $\overline{V^{fin}_q(\lambda)}$ is a module for $\bu$ and 
is isomorphic to $V^{fin}(\lambda)$.
It is also known \cite{L1} that any finite-dimensional $\bu_q^{fin}$-module 
$V_q$ is 
a direct sum of irreducible modules; we let $m_\mu(V_q)$ be the multiplicity 
with which $V_q^{fin}(\mu)$ occurs in $V_q$.

The type 1 irreducible finite--dimensional $\bu_q$-modules are parametrized by 
$n$-tuples of polynomials $\bpi_q =(\pi_1(u),\cdots, \pi_n(u))$, where the 
$\pi_r(u)$ have coefficients in $\bc(q)$ and constant term 1. Let us denote the 
corresponding module by $V_q(\bpi_q)$. Then, \cite{CP3}, there exists a unique 
(up to scalars) element $v_{\bpi_q}\in V_q(\bpi_q)$ satisfying 
\begin{equation}\label{qrel1}\bx_{k,r}^+.v_{\bpi_q} =0,\ \ K_i.v_{\bpi_q} 
=q^{{\text{deg}}\pi_i}v_{\bpi_q},\end{equation}
and 
\begin{equation}\label{qrel2} \bh_{i,k}. v_{\bpi_q} = d_{i,k}.v_{\bpi_q},\ \ 
(\bx_{i,k}^-)^{{\text{deg}}\pi_i+1}.v_{\bpi_q} =0,\end{equation}
where the $d_{i,k}$ are determined from the functional equation
\begin{equation*}{\text {exp}}\left(-\sum_{k\ge 0}\frac{d_{i,\pm 
k}u^k}{k}\right) = \pi_i^\pm(u),\end{equation*}
where $\pi_i^+(u) =\pi_i(u)$ and $\pi_i^-(u) = u^{\text{deg}\pi_i}\pi_i(u^{-1})/
\left.\left(u^{\text{deg}\pi_i}\pi_i(u^{-1}\right)\right|_{u=0}$.
We remark that these are in general {\bf not} the defining relations of 
$V_q(\bpi_q)$. Set,
\begin{equation*} V_\ba(\bpi_q) =\bu_\ba.v_{\bpi_q}.\end{equation*}

\begin{prop}{\label{afree}} Suppose that the $n$--tuple $\bpi_q=(\pi_1,\cdots \pi_n)$ is such that for all $j\in I$, $\pi_j(u)\in\ba[u]$. Regarded as an $\ba$--module $V_\ba(\bpi_q)$ is free of rank equal to $\text{dim}_{\bc(q)} V_q(\bpi_q)$.\end{prop}
\begin{pf}  In the simply--laced case, this was proved in \cite[Proposition  4.4]{CP5}.  The argument given there can be extended to include the case of $B_n$ and $C_n$ as follows. The crucial step is to prove that an element of the form
\begin{equation*}(\bx^-_{i_1,k_1})^{(s_1)}(\bx^-_{i_2,k_2})^{(s_2)}\cdots (\bx^-_{i_l,k_l})^{(s_l)}.v_{\bpi_q},\end{equation*}
can be rewritten as an $\ba$--linear combination of elements 
\begin{equation*}(\bx^-_{i'_1,k'_1})^{(s'_1)}(\bx^-_{i'_2,k'_2})^{(s'_2)}\cdots (\bx^-_{i'_l,k'_l})^{(s'_l)}.v_{\bpi_q},\ \ 0\le k'_j\le N(\eta)\end{equation*}
where $N(\eta)$ depends only on $\eta=\sum_js_j\alpha_{i_j}$ and $\bpi_q$.
The proof proceeds by an induction on $\text{ht}\ \eta$,  The case $\eta=s\alpha_{i}$ was done in \cite{CP5}. So we can assume that $s_{1}\ne 0$ and $s_{2}\ne 0$ and that $k_j\le N(\eta-s_1\alpha_{i_1})$ for all $2\le j\le l$. If $a_{i_1,i_2} =0$  the result is obvious. If $a_{i_1,i_2}=-1$ then the inductive step  is proved in \cite{CP5}. 

It remains to prove the inductive step when $a_{i_1,i_2}=-2$.   We assume $k_1\ge 0$, (the case $k_1<0$ is similar, see \cite{CP5}) and proceed by induction on $k_1$, with induction beginning at $k_1=N(\eta-s_1\alpha_1)$.  By Lemma \ref{aij}, we see that the elements  $(\gamma_{k_1, k_2})^{(r)}$ and $(\gamma'_{k_1,k_2})^{(t)}$  belong to  the $\bu_\ba$ subalgebra generated by  the elements $\{(\bx^-_{i,m})^{(s)}: i\in I, s\in\bz^+,  0\le m\le N(\eta-s_1\alpha_1)+2\}$. Now using Lemma \ref{aij1} and  the induction hypothesis we see that the element $(\bx^-_{i_1,k_1})^{(s_1)}(\bx^-_{i_2,k_2})^{(s_2)}\cdots (\bx^-_{i_l,k_l})^{(s_l)}.v_{\bpi_q}$ can be rewritten as a linear combination of similar elements but with the $k_j\le N(\eta-s_1\alpha_1) +2$ for all $j$ thus completing the inductive step. 

To complete the proof of the proposition, observe that 
since the module is finite--dimensional over $\bc(q)$, 
\begin{equation*} (\bx^-_{i_1,k_1})^{(s_1)}(\bx^-_{i_2,k_2})^{(s_2)}\cdots (\bx^-_{i_l,k_l})^{(s_l)}.v_{\bpi_q}=0\end{equation*}
for all $l>>0$ and for all  but finitely many values of $s_1,s_2,
\cdots s_l$. 
It follows now that, there exists an integer $N\ge 0$ such that the elements
\begin{equation*} (\bx^-_{i_1,k_1})^{(s_1)}(\bx^-_{i_2,k_2})^{(s_2)}\cdots (\bx^-_{i_l,k_l})^{(s_l)}.v_{\bpi_q},\ \ 0\le k_j<N\end{equation*} span $\bu_\ba. v_{\bpi_q}$. This means that $V_\ba(\bpi_q)$ is a finitely generated $\ba$--module and hence is a free $\ba$ module. Since these elements also clearly span  $V_q(\bpi_q)$ over $\bc(q)$, the proposition follows. \end{pf}
  Given, $\bpi=(\pi_1,\pi_2,\cdots ,\pi_n)$ such that $\pi_j(u)\in \ba[u]$ for all $j\in I$, set
\begin{equation}{\label{specialize}}\overline{V_q(\bpi_q)} = \bc_1\otimes_\ba 
V_\ba(\bpi_q).\end{equation}
Let $\overline{\bpi_q}$ be the $n$-tuple of polynomials with coefficients in 
$\bc$ 
obtained by setting $q=1$ in the components of $\bpi_q$. Then, 
$\overline{V_q(\bpi_q)}$ is a $\bu$-module  generated by $1\otimes v_{\bpi_q}$ 
and satisfying the relations in (\ref{qrel1}) and (\ref{qrel2}) with the 
generators $\bx_{i,k}^\pm$  etc.
being replaced by their classical analogues.
Further, if we write
\begin{equation*} V_q(\bpi_q)=\bigoplus_{\mu\in P^+} m_\mu(V_q(\bpi_q)) 
V_q^{fin}(\mu),\end{equation*}
as $\bu_q^{fin}$-modules, then 
\begin{equation*} \overline{V_q(\bpi_q)}=\bigoplus_{\mu\in P^+} 
m_\mu(V_q(\bpi_q)) V^{fin}(\mu),\end{equation*}
as $\bu$-modules.

From now on, we shall only be interested in the following case. Thus, for $i\in 
I$, $m\ge 0$, $a\in\bc^\times$, let 
$\bpi_q(i,m,a)$ be the $n$-tuple of polynomials given by 
\begin{align*} \pi_j(u)&=1,\ \ \ \ {\text{if $j\ne i$}},\\
\pi_i(u)&=(1-au)(1-aq^{-2}u)\cdots (1-aq^{-2m+2}u).\end{align*}
We denote the corresponding $\bu_q$-module by $V_q(i,m,a)$. 
In the case when $a=1$, we  set $ V_q(i,m,1) =V_q(i,m)$. For all $a\in\bc^\times$ we let $v_{i,m}$ denote the vector $v_{\bpi_q(i,m,a)}$.

Given any connected subset $J\subset I$, let $\bu_q^J$ be the quantized enveloping algebra of $L(\frak g_J)$, this clearly maps to the subalgebra of $\bu_q$ generated by the elements $\{\bx_{j,k}^\pm :j\in J, k\in\bz\}$.

\begin{lem} Let $J=\{i\}$,  $m\ge 0$. Then  $\bu_{J,q}.v_{i,m}\subset V_q(i,m)$ is an irreducible  $\bu_{J,q}$--module and  
 \begin{equation*}\bx_{i,k}^-. v_{i,m}= 
q^k\bx_{i,0}^-.v_{i,m}.\end{equation*}
In particular, 
\begin{equation*} {\text{dim}}_{\bc(q)}(V_q(i,m))_{m\lambda_i-\alpha_i} 
=1.\end{equation*} \end{lem}
\begin{pf}  It is easy to see that $\bu_{J,q}.w_{i,m}$ is an irreducible $\bu_{J_q}$--module. Further, a simple checking shows that the elements 
$\{\bx_{i,k}^-. v_{i,m}-
q^k\bx_i^-.v_{i,m}: k\in\bz\}$ generate a submodule of $V_q(i,m)$ not containing $v_{i,m}$ and hence must be zero.  
\end{pf}

In view of  ({\ref{specialize}) it follows from the preceding lemma, that
\begin{equation*} {\text{dim}}\overline{(V_q(i,m))}_{m\lambda_i-\alpha_i} 
=1.\end{equation*}
The next lemma is immediate.
\begin{lem} The $\bu$-module $\overline{V_q(i,m)}$ is a quotient of  $W(i,m)$.  
\hfill\qedsymbol\end{lem}
It now follows from Theorem \ref{main1}   that
\begin{lem}{\label{if}} For all $\mu\in P^+$ we have $m_\mu(V_q(i,m))\le 1$. Further,\begin{equation*} m_\mu(V_q(i,m))\ne 0\implies \mu\in P(i,m).\end{equation*}\hfill\qedsymbol\end{lem}

The main result of this paper is 
\begin{thm} {\label{main}}  Let $\mu\in P^+$. Then, $m_\mu(V_q(i,m))\le 1$ and $m_\mu(V_q(i,m))\ne 0$ if and only if $\mu\in 
P(i,m)$. \end{thm}

The following corollary is immediate.\begin{cor} 
For all $i\in I$ and $m\ge 0$, we have
\begin{equation*} W(i,m)\cong \overline{V_q(i,m)}.\end{equation*}
\hfill\qedsymbol\end{cor}

\noindent In view of Lemma \ref{if}, to prove Theorem \ref{main},  it suffices to prove
\begin{prop}\label{main2} Let $\mu\in P(i,m)$. Then, $m_\mu(V_q(i,m)) 
=1$.\end{prop} The rest of the section is devoted to proving this result.  
Observe that when $P(i,m)=\{m\lambda_i\}$ there is nothing to prove. This means that we can assume $\frak g$ is of type $B$, $C$ or $D$. 
 We  shall need 
the following result which is  a special case of a  theorem of \cite{K}, (see \cite{VV} for a different proof in the simply--laced case).  
\begin{prop}\label{cyclic} For all $i\in I$, $m\in\bz^+$,
the $\bu_q$--module $V_q(\lambda_i,1)\otimes 
V_q(\lambda_i, 
q^{-2})\otimes\cdots\otimes V_q(\lambda_i, q^{-2m+2})$ is generated by
$v_{i,1}\otimes v_{i,1}\otimes\cdots\otimes v_{i,1}$.\hfill\qedsymbol\end{prop}
Given two $n$-tuples of polynomials $\bpi_q$ and 
$\tilde{\bpi}_q$, let \begin{equation*}\bpi_q\ \tilde{\bpi}_q 
=(\pi_1\tilde{\pi}_1,\dots, \pi_n\tilde{\pi}_n).\end{equation*}

\begin{lem}\label{quotient} The assignment $v_{i,1}\otimes 
v_{i,1}\otimes\cdots\otimes v_{i,1}\mapsto v_{i,m}$ extends to a surjective homomorphism of  $\bu_q$--modules $phi_m^i: V_q(\lambda_i, 1)\otimes V_q(\lambda_i,q^{-2})\otimes\cdots\otimes 
V_q(\lambda_i,q^{-2m+2})\to V_q(i,m)$.  
\end{lem}
 \begin{pf} It was proved in \cite{CP3}, \cite{Da} that 
\begin{equation*}\Delta(\bh_{i,k})= \bh_{i,k}\otimes 1+1\otimes \bh_{i,k}\mod \bu\otimes\bu\bu(>)_+,\end{equation*}
where $\bu(>)$ is the subalgebra generated by $\bx_{j,l}^+$ for all $j\in I$ and $l\in\bz_+$, and $\bu(>)_+$ is the augmentation ideal. 
It is now easy to see using \eqref{qrel1} and \eqref{qrel2} that 
the action of \begin{equation*} \bh_{i,k}.v_{i,1}\otimes 
v_{i,1}\otimes\cdots\otimes v_{i,1}= (-1)^k q_i^{(\begin{matrix}\scriptstyle{m}\\ \scriptstyle{k}\end{matrix}t)}\left[\begin{matrix}m\\k\end{matrix}\right]_{q_i} v_{i,1}\otimes 
v_{i,1}\otimes\cdots\otimes v_{i,1},
\end{equation*}
and  
\begin{equation*} \bx_{i,k}^+.v_{i,1}\otimes 
v_{i,1}\otimes\cdots\otimes v_{i,1}=0,\end{equation*}
 for all $i\in I$ and $k\in\bz$. This proves the lemma. \end{pf}

To prove Proposition \ref{main2}, we proceed by induction on $m$. We first show that induction starts. 
 \begin{lem} {\label{fund}}\hfill
 
\begin{enumerate} 
\item[(i)] Assume $\frak g$ is of type $B_n$. If $i\ne n$, then as $\bu_q^{fin}$--modules we have  
\begin{align*} V_q(i,1)&\cong\bigoplus_{j=0}^{[i/2]} V_q^{fin}(\lambda_{i-2j}),\ \ i\ne n,\\ V_q(n,1)&\cong V_q^{fin}(\lambda_n),\\ V_q(n,2)& \cong V_q^{fin}(2\lambda_n)\bigoplus_{j=0}^{[n/2]} V_q^{fin}(\lambda_{n-2j}).\end{align*}
\item[(ii)]  Assume that $\frak g$ is of type $D_n$
and that $1\le i\le n-2$. Then, as $\bu_q^{fin}$-modules,
\begin{equation*} V_q(i,1) \cong\bigoplus_{j=0}^{[i/2]} V_q^{fin}(\lambda_{i-2j}),\end{equation*}
If $i=n-1,n$, then $V_q(i,1)\cong V_q^{fin}(\lambda_i)$.
\item[(iii)] If $\frak g$ is of type $C_n$, then
\begin{align*} V_q(i,1)&\cong V_q^{fin}(\lambda_i),\\
V_q(i,2)&\cong \oplus_{j=0}^iV_q^{fin}(2\lambda_j).\end{align*}
\end{enumerate}
 \end{lem}
\begin{pf} The case $m=1$ was proved in \cite{CP3}. Assume that $\frak g$ is of type $C_n$ and that $m=2$.
For $C_2$, the proposition was proved in  \cite{C}. Assume that we know the result for $C_{n-1}$. Take $J=\{2,\cdots ,n\}$.  By induction on $n$, we get
\begin{equation*} \bu_{J,q}.v_{i,m}=\bigoplus_{j=1}^i V_{J,q}^{fin}(2\lambda_j),\end{equation*}
(note that we regard $\lambda_j\in P^+$ as an element of $P_J^+$ by restriction). In other words, there exist vectors $0\ne w_j\in \left(\bu_{J,q}.w_{i,m}\right)_{2\lambda_j}$ for $1\le j\le i$ with 
\begin{equation*} E_{\alpha_r}.w_j =0\ \  \forall \ r\in J.\end{equation*}
Since $2\lambda_i-2\lambda_j\in\oplus_{i=2}^n\bz^+\alpha_i$, it follows that 
$E_{\alpha_1}.w_j=0$ as well. This proves that \begin{equation*} m_{2\lambda_j}(V(i,m))=1,\ \ \forall \ \ 1\le j\le i,\end{equation*}
 and 
hence it suffices to prove that the trivial representation occurs in $V_q(i,2)$. 
To prove this, let $K$ be the kernel of the map $\phi_2^i: V_q(i,1)\otimes V_q(i,1, q^{-2})\to V_q(i,2)$ defined in Lemma \ref{quotient}.   As $\bu_q^{fin}$--modules, we have 
$m_\mu(M) =1$ if $\mu=0$ or $\mu=2\lambda_1$. Let $w_0\in M$ be such that $E_{\alpha_j}.w_0= F_{\alpha_j}.w_0=0$ for all $j\in I$. Suppose that $w_0\in K$. Since $E_{\alpha_r}$ and $F_{\alpha_0}$ commute, we must have        $F_{\alpha_0}.w_0 =cw_1$  for some $0\ne c\in\bc^\times$. Since $w_1\notin K$ this means that $c=0$ and that $\bc.w_0$ is the trivial $\bu_q$--module.    This implies that the modules $V_q(i,1)$ and $V_q(i,1, q^2)$ are dual, but it is known, \cite{CP4}, that the dual of $V_q(i,1)$ is the module $V_q(i,1, q^d)$ where $d\ne 2$ is the Coxeter number of $C_n$. 
Hence $w_0\notin K$ and the multiplicity of the trivial module in $V_q(i,2)$ is one.  This proves the proposition for $C_n$.

The only remaining case is $B_n$ with $i=n$. But this is proved in the same way as for $C_n$. We omit the details.
\end{pf}

Given an $n$-tuple of polynomials $\bpi_q$,  and $J=\{2,\cdots ,n\}$, let $\bpi_{J,q} =(\pi_2,\cdots, 
\pi_n)$ and let $V_{J,q}(\bpi_{J,q})$ be the irreducible $\bu_{J,q}$-module  
associated to $\bpi_{J,q}$. Then, it is easy to see that  \begin{equation*}
\bu_{J,q}.v_{\bpi_q}\cong V_{J,q}(\bpi_{J_q}).\end{equation*} 
The next proposition was proved in \cite{CP4}, the proof is similar to the one given above for Lemma \ref{quotient}. 
\begin{prop}{\label{subalg}}\hfill 
The comultiplication $\Delta$ of $\bu_q$ induces   a $\bu_{J,q}$-module 
structure on $\bu_{J,q}. v_{\bpi_q}\otimes \bu_{J,q}. v_{\tilde{\bpi}_q}$.
Further, the natural map
\begin{equation*}\bu_{J,q}. v_{\bpi_q}\otimes \bu_{J,q}. v_{\tilde{\bpi}_q}\to  
V_{J,q}(\bpi_{J,q})\otimes V_{J,q}(\tilde{\bpi}_{J,q})\end{equation*}
is an isomorphism of $\bu_{J,q}$-modules (the right-hand side is regarded as  a 
$\bu_{J,q}$-module by using the comultiplication $\Delta_J$ of 
$\bu_{J,q}$)\hfill\qedsymbol
\end{prop}

It is clear from Lemma \ref{quotient} that there exists a $\bu_q$-module map 
$\phi^i_{m}: V_q(i,1)\otimes V_q((m-1)\lambda_i, q^{-2}) \to V_q(i,m)$ which 
maps $v_{i,1}\otimes v_{i,m-1}$ to $v_{i,m}$. For $J=\{2,3,\cdots ,n\}$, let  
$\phi^i_{J,m}$ denote the analogous map for $\bu_{J,q}$. From Proposition 
\ref{subalg}, we see that  the restriction of $\phi^i_{m}$ to 
$V_{J,q}(i,1)\otimes V_{J,q}((m-1)\lambda_i, q^{-2})$ is $\phi^{i-1}_{J,m}$.

In what follows we set $\phi_m=\phi^i_m$ and $\phi_{J,m}=\phi^{i}_{J,m}$ and we 
take $J=\{2,3,\cdots ,n\}$.

\noindent 
\noindent Proposition \ref{main2} follows from 
\begin{prop}{\label{phi}} Assume that $\frak g$ is of type $D_n$. Let $i\in I$, $i\ne 1, n-1, n$, and $m\ge 0$. For 
every $\mu\in P(i,m)$, there exist unique (up to scalars) non-zero elements 
$v^m_\mu\in V_q(i,m)_\mu$ 
with the following properties.

\begin{enumerate}

\item[$(i)_{m-1}$] If $\mu_1\in P_{i,1}$ and $\mu_2\in P_{i,m-1}$ are such that
$\mu_1+\mu_2=\mu$, then for some $c_{\mu_1,\mu_2}^\mu\in\bc^\times$,
\begin{equation*}\phi_m(v^1_{\mu_1}\otimes v^{m-1}_{\mu_2})= c_{\mu_1,\mu_2}^\mu 
v^m_\mu.\end{equation*}
\noindent \item[$(ii)_m$] For all $j\in I$,\begin{equation*} 
E_{\alpha_j}.v^m_\mu 
=0.\end{equation*} Further, if $\mu\in P(i,m)$ is such that $\mu+\lambda_2\in 
P(i,m)$, then 
\begin{equation*} F_{\alpha_0}. v^m_\mu=a_\mu 
v^m_{\mu+\lambda_2},\end{equation*}
for some $a_\mu\in\bc^\times$.
\end{enumerate}

Analogous statements hold for $B_n$ if $i\ne n$. If $i=n$ or if $\frak g$ is of type $C_n$, then we assume that  $m\ge 3$ and in $(i)_{m-1}$ that the  element $\mu_1\in P(i,2)$.
\end{prop}
\begin{pf} We begin by remarking that, if elements $v^m_\mu$ exist with the 
desired properties, then by Lemma \ref{if}, they are unique up 
to scalars. We shall only prove the proposition when $\frak g$ is of type $D_n$, the modifications in the other cases are clear.

Notice that by Lemma \ref{if} (since $E_{\alpha_r}$ and 
$F_{\alpha_0}$ commute for all $r\in I$), if $\mu\in P(i,m)$ is such that 
$\mu+\lambda_2\notin P(i,m)$, then \begin{equation*} F_{\alpha_0}.v_\mu^m 
=0.\end{equation*}
Also  observe  that if $\mu\in P(i,m)$, then $\mu+\lambda_2\in P(i,m)$ if 
and only if $\ell(\mu)<m$.  We shall use these facts throughout the proof with 
no further comment.

The statement $(i)_0$ is trivially true. For $(i)_1$ observe  that by 
Lemma \ref{fund}, we have non-zero vectors $v^1_\mu$ for $\mu\in P_{i,1}$ such 
that $ E_{\alpha_j}.v^1_\mu =0$ for all $j\in I$. If $i$ is even,
the only element $\mu\in P_{i,1}$ such that $\mu+\lambda_2\in P_{i,1}$ is 
$\mu=0$, and then we have
\begin{equation*} E_{\alpha_r}.v^1_0= F_{\alpha_r}.v^1_0 =0\ \ (r\in I). 
\end{equation*}
Thus, we have to prove that $F_{\alpha_0}. v^1_0\ne 0$. But this is 
clear, since
\begin{equation*} F_{\alpha_0}. v^1_0=0\implies E_{\alpha_0}. 
v^1_0=0,\end{equation*}
which would imply that $v^1_0$ generates a proper $\bu_q$-submodule of 
$V_q(i,1)$, 
contradicting the irreducibility of $V_q(i,1)$. If $i$ is odd, then 
$\mu+\lambda_2$ is not in $P_{i,1}$ for any $\mu\in P_{i,1}$ and hence the 
proposition is proved for $m=1$.

Assume from now on  that $(ii)_{m-1}$ and $(i)_{m-1}$ are known for $i$. We 
first prove that $(ii)_m$ and $(i)_m$ hold if $i$ is even. For $\mu\in P(i,m)$, 
let $\mu_1\in P_{i,1}$
and $\mu_2\in P_{i,m-1}$ be such that 
\begin{equation*}\mu=\mu_1+\mu_2.\end{equation*}
Set \begin{equation*} v^m_\mu=\phi_m(v_{\mu_1}^1\otimes 
v_{\mu_2}^{m-1}).\end{equation*}
Then, $v_\mu^m\ne 0$ since $(i)_{m-1}$ holds. Clearly,
\begin{equation*} E_{\alpha_j}.v^m_\mu= \phi_m(E_{\alpha_j}.(v^1_{\mu_1}\otimes 
v^{m-1}_{\mu_2})) =0.\end{equation*}
Suppose that $\mu+\lambda_2\in P(i,m)$, i.e., $\ell(mu)<m$. Then, either 
$\mu_1+\lambda_2\in P_{i,1}$ or $\mu_2+\lambda_2\in P_{i,m-1}$. For $j=1,2$, let 
$r_j=m-\ell(\mu_j)$. 
 Then, $F_{\alpha_0}^{r_1}.v^1_{\mu_1} =av^1_{\mu_1+r_1\lambda_2}$ and 
$F_{\alpha_0}^{r_2}.v^{m-1}_{\mu_2}=bv^{m-1}_{\mu_2+r_2\lambda_2}$ for some 
non-zero scalars $a,b\in\bc(q)$.

 Hence, 
\begin{equation*}F_{\alpha_0}^{r_1+r_2}.v^m_\mu
=\phi_m(F_{\alpha_0}^{r_1}v_{\mu_1}\otimes 
F_{\alpha_0}^{r_2}v_{\mu_2}).\end{equation*}
Since the right-hand side of the preceding equation is a non-zero scalar 
multiple of $v_{\mu+(r_1+r_2)\lambda_2}$, it follows that 
\begin{equation*} F_{\alpha_0}.v^m_{\mu}\ne 0.\end{equation*} 
Since $E_{\alpha_r}F_{\alpha_0}.v^m_\mu=0$ for all $r\in I$, it follows from 
Lemma \ref{if} that $F_{\alpha_0}.v_\mu= a_\mu v^m_{\mu+\lambda_2}$ for some 
non-zero scalar $a_\mu\in\bc(q)$. This shows that $(ii)_m$ holds when $i$ is 
even. To prove $(i)_m$, let $\mu_1\in P_{i,1}$ and $\mu_2\in P(i,m)$ and choose 
$r_1,r_2$ so that $\ell(\mu_1+r_1\lambda_2)=1$ and $\ell(\mu_2+r_2\lambda_2)=m$.
 Then,  \begin{equation*}
F_{\alpha_0}^{r_1}.v^1_{\mu_1}=v^1_{\lambda_2},\ \ 
F_{\alpha_0}^{r_2}.v^m_{\mu_2}=v^m_{\mu_2+r_2\lambda_2}.\end{equation*}
If $i=2$, we see that
\begin{equation*} F_{\alpha_0}^{r_1+r_2}.(v^1_{\mu_1}\otimes v^m_{\mu_2}) 
=v_{\lambda_2}^1\otimes v_{m\lambda_2}^m,
\end{equation*}
and hence that
\begin{equation*}
\phi_{m+1}( F_{\alpha_0}^{r_1+r_2}.(v^1_{\mu_1}\otimes v^m_{\mu_2})) 
=v^{m+1}_{(m+1)\lambda_2}.\end{equation*}
Clearly, this implies that $\phi_{m+1}(v^1_{\mu_1}\otimes v^m_{\mu_2})\ne 0$, 
and $(i)_m$ is proved when $i=2$. In particular, the theorem is proved for 
$n=4$. 

Assume that we know the proposition for $J=\{2,3,\cdots ,n\}$. 
Since $m\lambda_i-\mu_2-r_2\lambda_2\in Q_J^+$ and $\lambda-\lambda_2\in Q_J^+$, 
we see by the induction hypothesis on $n$ that
\begin{equation*}\phi_{m+1}(v^1_{\lambda_2}\otimes 
v^m_{\mu_2+r_2\lambda_2})=\phi_{J,m+1} (v^1_{\lambda_2}\otimes 
v^m_{\mu_2+r_2\lambda_2})\ne 0,\end{equation*}
i.e., that $\phi_{m+1}(F_{\alpha_0}^{r_1+r_2}.(v^1_{\mu_1}\otimes 
v^m_{\mu_2}))\ne 0$. This implies that $\phi_{m+1}(v_{\mu_1}\otimes 
v_{\mu_2})\ne 0$ and proves  that $(i)_m$ holds for $I$.

It remains to prove the result when $i$ is odd; recall that the proposition is 
known for $J$.
If $i$ is odd, then \begin{equation*}\mu\in P(i,m)\implies \ell(\mu)=m 
\implies\mu=m\lambda_i-\eta,\ \ (\eta\in Q_J^+).\end{equation*}
By the induction hypothesis, there exist elements $v_\mu\in 
\bu_{J,q}.v_{i,m}$ satisfying
\begin{equation*} E_{\alpha_j}.v^m_\mu =0\ \ {\text{for all}}\ j\in 
J.\end{equation*}
Clearly, $E_{\alpha_1}.v^m_\mu =0$, and this proves $(ii)_m$ since 
$\mu+\lambda_2$ is never in $P(i,m)$ if $i$ is odd. To see that $(i)_m$ holds, 
let $\mu_1\in P_{i,1}$ and $\mu_2\in P(i,m)$.
Then, $\mu_1\in\lambda_i-Q_J^+$ and $\mu_2\in m\lambda_i-Q_J^+$, and hence 
$v^1_{\mu_1}\in\bu_{J,q}v_{i,1}$ and  $v^m_{\mu_2}\in \bu_{J,q}.v_{i,m}$. Hence, 
\begin{equation*} \phi_{m+1}(v^1_{\mu_1}\otimes 
v^m_{\mu_2})=\phi_{J,m+1}(v^1_{\mu_1}\otimes v^m_{\mu_2})\ne 0, \end{equation*}
thus proving $(i)_m$ when $i$ is odd. The proof of the proposition is now 
complete.

\end{pf}

\section{The  Exceptional Algebras}  We summarize here the results that can be proved for the exceptional algebras, using the techniques and results of the previous sections.  Again we assume that the nodes are numbered as in \cite{B}.

\noindent $E_6$. Here $i\ne 4$. 

\begin{align*} V_q(i,m)& =V_q^{fin}(\lambda_i), \ \ i=1,6,\\
V_q(2,m)&\cong\bigoplus_{0\le r\le m} V_q^{fin}(r\lambda_2),\\
V_q(3,m) &\cong\bigoplus_{ r+s= m} V_q^{fin}(r\lambda_3+ s\lambda_6),\\
V_q(5,m) &\cong\bigoplus_{r+s= m} V_q^{fin}(r\lambda_5+ s\lambda_1).\end{align*}
\vskip 12pt

\noindent $E_7$. Here $i=1,2,6,7$.

\begin{align*} 
V_q(1,m)&\cong\bigoplus_{0\le r\le m}V_q^{fin}(r\lambda_1),\\
V_q(7,m)& \cong V_q^{fin}(\lambda_7), \\
V_q(2,m)&\cong\bigoplus_{r+s= m} V_q^{fin}(r\lambda_2+s\lambda_7),\\
V_q(6,m) &\cong \bigoplus_{0\le r+s\le  m}V_q^{fin}(r\lambda_6+ 
s\lambda_1).
\end{align*}

\medskip

\noindent $E_8$. Here $i=1,8$.
\begin{align*} V_q(1,m)&\cong\bigoplus_{0\le r\le  m} V_q^{fin}(r\lambda_8), \\
V_q(8,m) &\cong\bigoplus_{0\le r+s\le  m} V_q^{fin}(r\lambda_1+ s\lambda_8).\\
\end{align*}
w

\noindent $F_4$. 
\begin{align*} V_q(1,m)&\cong \bigoplus_{k=0}^m V_q^{fin}(s\lambda_1),\\
V_q(4,m)&\cong\bigoplus_{j=0}^k\oplus_{k=0}^{m/2} V_q^{fin}(j\lambda_1+(m-2k)\lambda_4).\end{align*}

\noindent $G_2$. 
\begin{equation*} V_q(1,m)\cong\bigoplus_{k=0}^m V_q^{fin}(k\lambda_1).\end{equation*}


\begin{thebibliography}{BCPQR}

\bibitem[AK]{AK} T.~Akasaka and M.~Kashiwara, Finite-dimensional 
representations of quantum affine algebras, {\it Publ. Res. Inst. Math. Sci.} 
{\bf 33} (1997), no. 5, 839-867.


\bibitem[B]{B} J.~Beck,
Braid group action and quantum affine algebras, {\it Commun. Math. Phys.} {\bf
165} (1994), 555-568.

\bibitem[BCP]{BCP} J.~Beck, V.~Chari and A.~Pressley,
An algebraic characterization of the affine canonical basis,
{\it Duke Math. J.} {\bf 99} (1999), no. 3, 455-487.

\bibitem[B]{Bo} N.~Bourbaki,  Groupes et alg\`ebres de Lie, Chapitres
4,5,6,  Hermann, Paris (1968). 




\bibitem[C]{C} V.~Chari, Minimal affinizations of representations of affine Lie 
algebras: the rank 2 case, {\it Publ. Res. Inst. Math. Sci.} {\bf 31} (1995), 
no. 5, 893-911.

\bibitem[CP1]{CP1} V.~Chari and A.~Pressley,  Quantum affine algebras, {\it 
Commun. Math. Phys.} {\bf 142} (1991), 261-283.

\bibitem[CP2]{CP2} V.~ Chari and A.~Pressley, {\it A Guide to Quantum 
Groups}, Cambridge University Press, Cambridge, 1994.

\bibitem[CP3]{CP3} V.~Chari and A.~Pressley, Quantum affine algebras and their 
representations, in Representations of Groups, (Banff, AB, 1994), 59-78, {\it 
CMS Conf. Proc.} {\bf 16}, AMS, Providence, RI, 1995.

\bibitem[CP4]{CP4} V.~Chari and A.~Pressley,  Minimal affinizations of 
representations of quantum groups: the simply-laced case, {\it J. Alg.} {\bf 
184} (1996), no. 1, 1-30.

\bibitem[CP5]{CP5}  V.~Chari and A.~Pressley, Weyl modules for classical and 
quantum affine algebras, preprint, q-alg/0004174.

 \bibitem[Da]{Da} I. Damiani, {\em La $R$-matrice pour les alg\`ebres
    quantiques de type affine non tordu}, preprint.


 



\bibitem[Dr1]{Dr1} V.G.~Drinfeld, Hopf Algebras and the quantum Yang-Baxter 
equation, {\it Sov. Math. Dokl.}  {\bf 32} (1985), 254-258.

\bibitem[Dr2]{Dr2} V.G.~Drinfeld,
A new realization of Yangians and  quantum
affine algebras.
{\it Soviet Math. Dokl.} {\bf 36} (1988), 212-216.



\bibitem[GV]{GV} V.~Ginzburg and E.~Vasserot, Langlands reciprocity for affine 
quantum groups of type $A_n$, {\it Int. Math. Res. Not.} {\bf 3} (1993), 67-85.

\bibitem[J]{J} N.~Jing,
On Drinfeld realization of quantum affine algebras. The
Monster and Lie algebras (Columbus, OH, 1996), pp. 195-206,
{\it Ohio State Univ. Math. Res. Inst.
Publ.} {\bf 7}, de Gruyter, Berlin, 1998.

\bibitem[FM]{FM} E.~Frenkel and E.~Mukhin, Combinatorics of $q$-characters of 
finite-dimensional representations of quantum affine algebras, preprint, 
math.qa/9911112.

\bibitem[FR]{FR} E.~Frenkel and N.~Reshetikhin,  The $q$-characters of 
representations of quantum affine algebras and deformations of $W$-algebras,
{\it Contemp. Math} {\bf 248} (1999).

\bibitem[HKOTY]{HKOTY} G.~Hatayama, A.~Kuniba, M.~Okado, T.~Takagi, Y.~Yamada,  Remarks on the Fermionic Formula, {\it Contemp. Math} {\bf 248}.


\bibitem[K]{K} M.~Kashiwara, On level zero representations of quantized affine algberas, qa/0010293/ 
\bibitem[KR]{KR} A.N.~Kirillov and N.~Reshetikhin, Representations of Yangians 
and multiplicities of ocurrence of the irreducible components of the tensor 
product of simplie Lie algebras, {\it J. Sov. Math.} {\bf 52} (1981), no. 5, 
393-403.

\bibitem[KSS]{KSS} A.N.~Kirillov, A. Schilling and M. Shimozono, A bijection between Littelwood--Richardson tableaux and rigged congigurations, math.CO/9901037.


\bibitem[Kl]{Kl} M.~Kleber, Combinatorial structure of finite-dimensional 
representations of Yangians: the simply-laced case, {\it Int. Math. Res. Not.}
{\bf 7} (1997), no. 4, 187-201.





\bibitem[L1]{L1} G.~Lusztig, Quantum deformations of certain simple modules over 
enveloping algebras, {\it Adv. Math.} {\bf 70} (1988), 237-249.
\bibitem[L2]{L3} G.~Lusztig, Quantum groups at roots of 1, Geom. Ded. {\bf{35}} (1990), 89--114.

\bibitem[L3]{L2} G.~Lusztig, {\em Introduction to quantum groups}, Progress in
Mathematics {\bf 110}, Birkh\"auser,
Boston, 1993.
\bibitem[VV]{VV} M.~Varagnolo and E.~ Vasserot, Standard modules for quantum affine algberas, preprint.






\end{thebibliography}
\end{document}